\documentclass[12pt,leqno]{amsart}
\usepackage[dvips]{graphicx}
\usepackage{amsfonts}
\usepackage{amssymb}
\usepackage{amsmath}
\usepackage{amscd}
\usepackage{graphicx}
\def\DATE{\today}
\newtheorem{theorem}{Theorem}
\newtheorem{definition}[theorem]{Definition}
\newtheorem{corollary}[theorem]{Corollary}

\newtheorem{proposition}[theorem]{Proposition}

\newcommand\R{\mathbb{R}}

\newcommand\g{\mathfrak{g}}

\newcommand\K{\mathbb{K}}

\newcommand\pf{\noindent{\it Proof. }}
\newcommand\lr{\left\{ \begin{array}{l}}
\def\ds{\displaystyle}

\title{Weakly associative and symmetric Leibniz algebras }
\author{Elisabeth Remm}
\date{Chevat 16, 5775}
\address{Universit\'{e} de Haute Alsace, LMIA, 6 rue des Fr\`{e}res Lumi\`{e}re, 68093 Mulhouse}
\email{elisabeth.remm@uha.fr}

\begin{document}

\maketitle

\begin{abstract}
We study a special class of weakly associative algebras:  the symmetric Leibniz algebras. We describe the structure of the commutative and skew symmetric algebras associated with the polarization-depolarization principle. We also give a structure theorem for the symmetric Leibniz algebras and we study formal deformations in the context of deformation quantization.

\end{abstract}

\medskip

\noindent{\bf Introduction}

The notion of weakly associative algebra was introduced in \cite{R-Wass} in order  of broadening the notion of  deformation quantization of Poisson algebras. In fact, a deformation quantization of a Poisson algebra $(A,\bullet,[\; ,\, ])$ is given by a formal
deformation of the associative algebra $(A,\bullet)$ whose linear term (the term of degree $1$ in the series expansion of the formal deformation) determines the Lie bracket $[\; ,\, ]$. Recall that $(A,\bullet, [\;,\, ])$ is a Poisson algebra if
\begin{enumerate}
\item $(A,\bullet)$ is a associative commutative algebra,
\item $(A,[\;,\,])$ is a Lie algebra,
\item these two multiplications on $A$ are linked by the Leibniz Identity
$$[X\bullet Y,Z]=X\bullet [Y,Z]+[X,Z] \bullet Y$$
for any $X,Y,Z \in A.$
\end{enumerate}
In this context, a formal deformation $(A_t,\bullet_t)$ of $(A,\bullet)$ is a formal associative algebra and the set of all these deformations are parametrized by the Hochschild cohomology of $(A,\bullet).$
It has been shown in \cite{R-Wass} that a Poisson algebra $(A,\bullet,[\; ,\,])$ can be obtained  from a formal deformation of the commutative associative algebra $(A,\bullet)$ but in a larger category of algebras namely weakly associative. This class of algebras contains in particular the Lie algebras, the associative algebras and the commutative algebras. We show that the {\it symmetric Leibniz algebras}, that is algebras $(A,\ast)$ satisfying the pair of identities
$$
\left\{
\begin{array}{l}
X\ast(Y\ast Z) = (X\ast Y)\ast Z + Y\ast(X\ast Z),\\
 (Y\ast Z)\ast X = (Y\ast X)\ast Z + Y\ast (Z\ast X)
\end{array}
\right.
$$
for any $X,Y,Z \in A$, are also {\it weakly associative}. We use this property to describe these algebras which leads to a structure theorem. As an application we find  in particular  some results of \cite{A.B}.

In the first section we recall definition and properties of weakly associative algebras, the polarization-depolarization principle which permits to view these algebras as nonassociative Poisson algebras. We then in Section 2 study the properties of symmetric Leibniz firstly as one multiplication nonassociative algebra: the are weakly associative algebras. But we also study symmetric Leibniz algebras viewed as Poisson algebra after polarization and show that the associative commutative product $\bullet$ is 2-step nilpotent and the Lie bracket $[\; , \, ]$ has the property that
$[X,Y\bullet Z]=0$ and $X \bullet [Y, Z]=0$ for any vectors $X,Y,Z$. Then in the third section we give a classification of symmetric Leibniz algebras for the small dimensional cases.

\section{Weakly associative algebras}
Let $\K$ be a field of characteristic $0$ and $(A,\ast)$ a $\K$-nonassociative algebra. Recall  that an algebra $(A,\ast)$ is a nonassociative algebra  if it is a $\K$-vector space   with a bilinear binary multiplication $X \ast Y$ which may or not be associative.
The associator of the algebra $(A,\ast)$ is the trilinear map  $\mathcal{A}_\ast$ defined by
$$\mathcal{A}_\ast (X,Y,Z)=X \ast (Y \ast Z)- (X \ast Y )\ast Z$$
for any $X,Y,Z \in A.$
\subsection{Definition and examples of weakly-associative algebras}
\begin{definition}\cite{R-Wass}
A $\K$-algebra $(A,\ast)$ is called weakly associative if its associator $\mathcal{A}_\ast$ satisfies
$$\mathcal{A}_\ast (X,Y,Z)+\mathcal{A}_\ast (Y,Z,X)-\mathcal{A}_\ast (Y,X,Z)=0$$
for any $X,Y,Z \in A.$
\end{definition}
To simplify, we shall denote by $\mathcal{WA}_\ast (X,Y,Z)$ the trilinear map
$$\mathcal{WA}_\ast (X,Y,Z)=\mathcal{A}_\ast (X,Y,Z)+\mathcal{A}_\ast (Y,Z,X)-\mathcal{A}_\ast (Y,X,Z)$$
and a $\K$-algebra $(A,\ast)$ is weakly associative  if and only if
$$\mathcal{WA}_\ast (X,Y,Z)=0$$
for any $X,Y,Z \in A.$
\medskip

\noindent{\bf Examples.}
\begin{enumerate}
  \item Any associative  algebra is weakly associative.
  \item Any abelian (i.e. commutative) algebra is weakly associative.
  \item Any Lie algebra is weakly associative. In fact, if $(\g, \{,\})$ is a $\K$-Lie algebra, the associator of the Lie bracket satisfies
  $$\mathcal{A}_{[\, , \,  ]} (X,Y,Z)=[X,[Y,Z]\,]-[\,  [X,Y],Z]$$
  and from the Jacobi identity
  $$\mathcal{A}_{[\, , \,  ]}  (X,Y,Z)=[\, [Z,X],Y].$$
  Then
  $$\mathcal{WA}_{[\, ,\, ]}  (X,Y,Z)=[\, [Z,X],Y]+[\, [X,Y],Z]+[\, [Y,Z],X]$$
  and, from the Jacobi identity
   $$\mathcal{WA}_{[\, ,\, ]}  (X,Y,Z)=0$$
	\end{enumerate}
In the next section, we will show another important example of weakly associative algebras given by the class of {\it symmetric Leibniz algebras}.  As said in the introduction, the notion of weakly associative algebra has been defined in relation with the classical notion of Poisson algebra. This relation is based on the concept of polarization-depolarization presented in \cite{MRPoisson} which presents Poisson algebras in term of algebras with one nonassociative multiplication so in term of nonassociative algebras. 

We can also recall a characterization of weakly associative algebra

\begin{proposition}\label{der}\cite{R-Wass}
Let $(A,\ast)$ be a nonassociative $\K$-algebra. Then for any $X \in A$, the endomorphism
$L_X-R_X$ defined by
$$(L_X-R_X)(Y)=X \ast Y - Y \ast X$$
is a derivation of $(A,\ast)$ if and only if $(A,\ast)$ is weakly associative.
\end{proposition} 

\subsection{Polarization-Depolarization of $\K$-algebras} The polarization technique consists in representing a given one-operation $\K$-algebra $(A, \ast)$ without particular symmetry as an algebra with two operations, one commutative
and the other skew symmetric. Explicitly, in the following, we will decompose the multiplication $\ast$ of the $\K$-algebra $(A, \ast)$ using:
\begin{enumerate}
  \item $X\bullet Y= \frac{1}{2}(X\ast Y + Y \ast X)$ its symmetric part,
  \item $[X,Y]= \frac{1}{2}(X\ast Y - Y \ast X)$ its skewsymmetric part,
\end{enumerate}
for $X,Y \in A$. 
The triplet $(A, \bullet , [ \, \,  , \, ])$  will be referred to as the {\it polarization} of $(A, \ast)$. Conversely, starting with an algebra $(A, \bullet , [ \, \,  , \, ])$ with a  commutative product $\bullet$ and a skew symmetric multiplication $[ \, \,  , \, ]$, we obtain an algebra $(A, \ast)$ with only one nonassociative multiplication defined by $X \ast Y= X\bullet Y + [X ,Y]$. The algebra $(A, \ast)$ is called the {\it depolarization} of the algebra $(A, \bullet , [ \, \,  , \, ])$.
The polarization-depolarization technique of \cite{MRPoisson} then makes a link between Poisson algebras and some nonassociative algebras and permits to present Poisson algebras as a nonassociative algebras $(A,\ast)$ satisfying 
$$3(x\ast y)\ast z-3x\ast (y\ast z) = (x\ast z)\star y + (y\ast z)\ast x − (y\ast x)\ast z − (z\ast x)\ast y.$$
Another applications of this technique are found in \cite{Ben}.

The polarization also permits to study weakly associative algebras as nonassociative Poisson algebras.

\begin{definition}
Let  $(A,\bullet,[\, ,\, ])$ be  a triple where $A$ is a $\K$-vector space with two multiplications  $\bullet$ and $  [\, ,\, ]$. We say that $(A,\bullet,[\, ,\, ])$ is a nonassociative Poisson algebra if
  \begin{enumerate}
  \item $(A,\bullet)$ is a  nonassociative commutative algebra,
  \item $(A,[\, ,\, ] )$ is a Lie algebra,
	\item the Leibniz identity between $[\, ,\, ]$ and $\bullet$ is satisfied, that is
$$[X\bullet Y,Z]=X\bullet[Y,Z]+[X,Z]\bullet Y$$
fo all  $X,Y,Z \in A.$
\end{enumerate}
\end{definition}
\noindent{\bf Remark.}  As the terminology suggests, nonassociative Poisson algebras generalize Poisson algebras by relaxing the
associative condition on the underlying commutative algebra.

The link between Poisson algebras and weakly associative algebra is summarized in the following result:
\begin{theorem}\cite{R-Wass}
Let  $(A,\bullet,[ \; , \, ])$ be a nonassociative Poisson algebra and consider on $A$ the third multiplication $$X \ast Y= X\bullet Y + [X,Y].$$
Then the algebra $(A, \ast)$, that is its depolarization, is weakly associative.
Conversely, 
if  $(A,\ast)$ is a weakly associative algebra, then its polarization  $A(\bullet,[\; ,\, ])$ is a nonassociative Poisson algebra.
\end{theorem}
\noindent{\bf Remark.}
 Let  $( A,\ast)$ a Lie-admissible algebra and $(A,\bullet, [ \; ,  \, ])$ its polarization. Then the following relations are equivalent
 \begin{enumerate}
  \item Leibniz($[ \; ,\, ],\ast$) : $[X\ast Y,Z]-X \ast [Y,Z]-[X,Z] \ast Y=0,$
  \item  Leibniz($[ \; ,\, ],\bullet$) : $[X\bullet Y,Z]-X\bullet[Y,Z]-[X,Z]\bullet Y=0,$
\end{enumerate}
for any $X,Y,Z \in A.$

Remark also that in order to obtain the classical notion of Poisson algebra, the multiplication $\bullet$ has to be associative. We have
$$\begin{array}{lll}
4\mathcal{A}_\bullet(X,Y,Z)&=& X \ast (Y \ast Z + Z \ast Y)+(Y \ast Z + Z \ast Y) \ast X-(X\ast Y + Y \ast X) \ast Z  \\
&&-Z \ast (X \ast Y + Y \ast X) \\
\end{array}$$If we denote by $\mathcal{B}_\ast$ the trilinear map associated to the multiplication $\ast$ defined by
$$\mathcal{B}_\ast(X,Y,Z)=X\ast(Y \ Z)+(Y \ast Z) \ast X$$
then
$$4\mathcal{A}_\bullet(X,Y,Z)= \mathcal{B}_\ast(X,Y,Z)+\mathcal{B}_\ast(X,Z,Y)-\mathcal{B}_\ast(Z,X,Y)-\mathcal{B}_\ast(Z,Y,X).$$
Then a sufficient condition for the associativity of $\bullet$ is that $\mathcal{B}_\ast(X,Y,Z)=0$.

\section{Symmetric Leibniz algebras}
A symmetric Leibniz algebra is an algebra $(A, \ast)$ such that for any $X,Y ,Z \in
A$, we have
$$
\left\{
\begin{array}{l}
X\ast(Y\ast Z) = (X\ast Y)\ast Z + Y\ast(X\ast Z),\\
 (Y\ast Z)\ast X = (Y\ast X)\ast Z + Y\ast (Z\ast X).
\end{array}
\right.
$$
If we denote by $\mathcal{A}_\ast$ the associator of the multiplication $ \ast$, that is
$$\mathcal{A}_\ast (X,Y,Z)=X \ast (Y \ast Z)- (X \ast Y ) \ast Z$$
for any $X,Y,Z \in A$, the first identity corresponds to
$$\mathcal{A}_\ast (X,Y,Z)=Y \ast (X \ast Z)$$
and the second to
$$\mathcal{A}_\ast (Y,Z,X)=-(Y\ast X)\ast Z.$$
We deduce that  $(A, \ast)$ is a symmetric Leibniz algebra if and only if
$$
\left\{
\begin{array}{l}
\mathcal{A}_\ast (X,Y,Z)=Y \ast (X \ast Z)\\
 \mathcal{A}_\ast (X,Y,Z)=-(X\ast Z)\ast Y
\end{array}
\right.
$$
or equivalently
$$
\left\{
\begin{array}{l}
\mathcal{A}_\ast (X,Y,Z)=Y \ast (X \ast Z)\\
\mathcal{B}_\ast(X,Y,Z)=0.
\end{array}
\right.
$$
In particular we deduce
$$\mathcal{A}_\ast (X,Y,Z)+\mathcal{A}_\ast (Y,Z,X)=Y \ast (X \ast Z)
-(Y\ast X)\ast Z=\mathcal{A}_\ast(Y,X,Z)$$
that is
$$\mathcal{A}_\ast (X,Y,Z)+\mathcal{A}_\ast (Y,Z,X)-\mathcal{A}_\ast(Y,X,Z)$$
and $A$ is a weakly associative algebra.
\begin{proposition}
Any symmetric Leibniz algebra is weakly associative.
\end{proposition}

Since any symmetric Leibniz algebra $(A,\ast)$ is weakly associative, it is also Lie admissible \cite{GRLieadm} and its polarized algebra $(A,\bullet,[\; , \, ])$ is a nonassociative Poisson algebra. But we also have that $\mathcal{B}_\ast(X,Y,Z)=0$ for any $X,Y,Z \in A$ so we obtain the following result:

\medskip

\begin{proposition}
Let $(A,\ast)$ be a symmetric Leibniz algebra. Then it is a weakly associative algebra and the commutative multiplication $\bullet$ of its polarization $A(\bullet,[\;,\,])$ is also associative. Then the polarization  $A(\bullet,[\;,\, ])$ is a  Poisson algebra.
\end{proposition}
So a Poisson algebra can be naturally obtained from each symmetric Leibniz algebra by polarization. We find a result of \cite{A.B}.

\medskip

\noindent Remark. In \cite{GRnonass}, we have studied some classes of nonassociative algebras in terms of action of the symmetric group $\Sigma_3$. The class of symmetric Leibniz algebras are a part of these, more precisely are $v-w$-algebras where $v,w $ are in the algebra group $\K[\Sigma_3]$. 
\section{Structure of symmetric Leibniz algebras}
Recall that, if $(A, \ast)$ is a symmetric Leibniz algebra, we have in particular the relation
$$\mathcal{B}_\ast(Y,X,Z)=Y\ast (X \ast Z)+(X \ast Z) \ast Y = 0.$$
This is equivalent, by polarization, to
$$\begin{array}{l}
[Y,[X,Z]]+[Y,X \bullet Z]+Y \bullet [X,Z]+Y \bullet (X \bullet Z)+[[X,Z],Y]+[X \bullet Z,Y]+[X,Z] \bullet Y\\
 +(X \bullet Z) \bullet Y=0
\end{array}
$$
which is equivalent to
$$
Y \bullet [X,Z]+Y \bullet (X \bullet Z)+[X,Z] \bullet Y
 +(X \bullet Z) \bullet Y=0
$$
or
$$Y \bullet [X,Z]+Y \bullet (X \bullet Z)=0.$$
We deduce that we have also
$$Y \bullet [Z,X]+Y \bullet (Z \bullet X)=0.$$
If we add these two identities we obtain:
\begin{proposition}
Let $(A,\ast)$ be a symmetric Leibniz algebra and $(A,\bullet,[\; ,\, ])$ its polarization. Then $\bullet$ is an associative commutative multiplication which satisfies
$$X \bullet (Y \bullet Z)=0$$
for any $X,Y,Z \in A$.
\end{proposition}
\begin{corollary}
The commutative associative algebra $(A,\bullet)$ is two-step nilpotent.
\end{corollary}
If we denote by $A^2_\bullet$ the subalgebra of $(A,\bullet)$ generated by the products $X\bullet Y$, $X,Y \in A$ and by $A^1_\bullet$ a vectorial subspace of $A$ isomorphic to $A/A^2_\bullet$ then we have the grading
$$A= A^1_\bullet \oplus A^2_\bullet$$
with $$A^1_\bullet \bullet A^1_\bullet  = A^2_\bullet, \ A^1_\bullet \bullet A^2_\bullet = A^2_\bullet \bullet A^2_\bullet =0.$$

\medskip

Now lets look at the properties of the Lie bracket associated with the polarization process of the symmetric Leibniz multiplication. We have seen that 
$Y \bullet [X,Z]+Y \bullet (X \bullet Z)=0$ for any $X,Y,Z \in A$. Since 
$Y \bullet (X \bullet Z)=0$, we deduce
\begin{proposition}
Let $(A, \ast)$ be a symmetric Leibniz algebra and $(A, \bullet, [\; , \, ])$ its polarization.
The derived Lie subalgebra of the Lie algebra $(A,[\; ,\, ])$, denoted by 
$[A,A],$ is contained in the center of the associative algebra $(A,\bullet)$.
\end{proposition}
Since for any $X,Y,Z \in A$ we have $X \bullet (Y \bullet Z)=0$ and $X \bullet [Y,Z]=0$, the relation $\mathcal{A}_\ast(X,Y,Z)-Y\ast (X \ast Z)=0$ is equivalent to
$$[X,Y \bullet Z] -[X \bullet Y,Z]-[Y,X\bullet Z]=0.$$
Then we have also
$$[X,Z \bullet Y] -[X \bullet Z,Y]-[Z,X\bullet Y]=0.$$
Since $\bullet$ is commutative, by adding these relations we obtain that
$$[X,Y \bullet Z]=0$$
for any $X,Y,Z \in A$.
\begin{theorem}
Let $(A,\ast)$ be a $\K$-algebra and $(A,\bullet,[\;,\,])$ its polarized algebra. Then $(A,\ast)$ is a symmetric Leibniz algebra if and only if the following conditions are satisfied:
\begin{enumerate}
\item $X \bullet (Y \bullet Z)=0$ for any $X,Y,Z \in A$

(implying that $(A,\bullet)$ is a commutative associative $2$-step nilpotent algebra).
\item $[\; ,\,]$ is a Lie bracket with the property that
$$[X,Y \bullet Z]=0 \ { and} \ X \bullet [Y,Z]=0$$
for any $X,Y,Z \in A$.
\end{enumerate}
\end{theorem}
\pf We have to prove only the converse. If $\bullet$ and $[\; ,\,]$ satisfy
$$X \bullet (Y \bullet Z)=0,\ \ [X,Y \bullet Z]=0, \ \ X \bullet [Y,Z]=0$$
then
$$
\begin{array}{lll}
2(\mathcal{A}_\ast(X,Y,Z)-Y\ast (X \ast Z))&=& X \bullet (Y \bullet Z)-(X \bullet Y )\bullet Z+X \bullet [Y,Z]+[X,Y]\bullet Z\\
&&+[X,Y\bullet Z]-[X\bullet Y,Z]+[X,[Y,Z]]+[Z,[X,Y]]\\
&& -Y\bullet (X \bullet Z)-Y \bullet [X ,Z]-[Y,X \bullet Z]+[Y,[Z,X]]\\
&=&0.
\end{array}
$$
The second identity of symmetric Leibniz algebras is shown in the same way.
\medskip

\noindent{\bf Consequence.} Since any symmetric Leibniz algebra $(A,\ast)$ is weakly associative satisfying  $\mathcal{B}_\ast(X,Y,Z)=0$ for any $X,Y,Z \in A$, its polarized algebra has a Poisson algebra structure. From the previous theorem, the Leibniz identity
$$[X\bullet Y, Z]-X\bullet [Y,Z]-[X,Z]\bullet Y=0$$
is trivially satisfied because each term is null. Let us also remark that $(A,\ast,[\; ,\, ])$ satisfies also Leibniz($[\; , \, ],\ast$) is also satisfied, that is 
$$[X\ast Y,Z]-X \ast [Y,Z]-[X,Z] \ast Y=0.$$

\section{Finite dimensional case}
Let $(A, \ast)$ a symmetric Leibniz $\mathbb{K}$-algebra and $(A,\bullet,[\; ,\,])$ its polarized Poisson algebra. In the finite dimensional case, we can find a basis $\{u_1,\cdots,u_r,v_1,\cdots,v_p,w_1,\cdots,w_q\}$ of $A$ such that $\{w_1,\cdots,w_q\}$ is a basis of $A^2_\bullet$ and $\{v_1,\cdots,v_p,w_1,\cdots,w_q\}$  a basis of the center of the associative algebra $(A,\bullet)$ (recall that $A^2_\bullet$ is contained in the center of $(A,\bullet)$). The nontrivial structure constants of $(A,\bullet)$ are given by
$$u_i \bullet u_j=\sum _{k=1}^q A_{i,j}^k w_k$$
with $A_{i,j}^k=A_{j,i}^k$.

Since $A^2_\bullet$ is contained in the center of $(A,[\; , \, ])$ and $A^2_{[\, ,\, ]}$ is contained in the center of $(A,\bullet)$, where $A^2_{[\, ,\, ]}$ is the derived Lie subalgebra of $(A,[\; ,\, ])$, we have
$$\begin{array}{l}
[u_i,w_j]=0, \ 1 \leq i \leq r, \ 1 \leq j \leq q,\\
\lbrack v_i,w_j]=0, \ 1 \leq i \leq p, \ 1 \leq j \leq q\\
\end{array}
$$
and
$$\begin{array}{l}
\medskip
[u_i,u_j]=\ds\sum_{k=1}^p C_{i,j}^kv_k+\sum_{l=1}^q D_{i,j}^lw_l,\\
\medskip
\lbrack u_i,v_j]=\ds\sum_{k=1}^p E_{i,j}^kv_k+\sum_{l=1}^q F_{i,j}^lw_l,\\
\medskip
\lbrack v_i,v_j]=\ds\sum_{k=1}^p G_{i,j}^kv_k+\sum_{l=1}^q H_{i,j}^lw_l\\
\end{array}
$$
with the Jacobi polynomial identities.  

The classifications in small dimensions of the commutative associative $2$-step nilpotent  algebras are already established \cite{dG, Rh}. Recall these results:
\begin{enumerate}
\item Dimension $2$
\begin{enumerate}
\item $A^2_\bullet=0$, that is $A=\K\{v_1,v_2\}$
$$v_1 \bullet v_2 =0.$$
\item $A=\K\{u_1,w_1\}$ and 
$$u_1 \bullet u_1 = w_1.$$
\end{enumerate}

\item Dimension $3$ (From now on, we don't write the decomposable algebras).
\begin{enumerate}
\item  $A=\K\{u_1,u_2,w_1 \}$
$$u_1 \bullet u_1 = {w_1}, \ u_2 \bullet u_2 = w_1$$
\item $\K=\R$, $A=\K\{u_1,u_2,w_1 \}$
$$u_1 \bullet u_1 = {w_1}, \ u_2 \bullet u_2 = -w_1$$
\item $A=\K\{u_1,u_2,w_1 \}$
$$u_1 \bullet u_2 =  u_2 \bullet u_1 = w_1.$$
\end{enumerate}
\item Dimension $4$ ($\K$ is algebraically closed)
\begin{enumerate}
\item  $A=\K\{u_1,u_2,u_3 \} \oplus \K\{w_1\}$ and
$$u_1 \bullet u_1 = w_1, \ u_2 \bullet u_2 =w_1, \ u_3 \bullet u_3 = w_1$$

\item $A=\K\{u_1,u_2 \} \oplus \K\{w_1,w_2\}$ and
$$u_1 \bullet u_1 = w_1, \ \ u_1 \bullet u_2=u_2 \bullet u_1=w_2$$

\item $A=\K\{u_1,u_2 \} \oplus \K\{w_1,w_2\}$ and
$$u_1 \bullet u_1 = w_1, \ \ u_2 \bullet u_2 =  w_1, \ \ u_1 \bullet u_2=u_2 \bullet u_1=w_2$$

\end{enumerate}
\end{enumerate}

\medskip

Using this list we describe all the small dimensional symmetric Leibniz algebras.  

\noindent{\bf 1. $\dim A=2.$}
\begin{enumerate}
\item $A^2_\bullet=0$, that is $A=\K\{v_1,v_2\}$ and 
$v_1 \bullet v_2 =0.$
In this case $[v_1,v_2]=av_1+bv_2$ and we obtain two classes of symmetric Leibniz algebras
$$
\left\{
\begin{array}{l}
v_i\ast v_j = 0, \ \forall i,j \in \{i,j\}
\end{array}
\right.
 \ \ \ \ 
\left\{
\begin{array}{l}
v_1 \ast v_2=-v_2 \ast v_1 =v_2.
\end{array}
\right.
$$
\item $A=\K\{u_1,w_1\}$ and 
$u_1 \bullet u_1 = w_1.$
Its center is $1$-dimensional and coincides with $A_\bullet ^2=\K\{w_1\}$.  We deduce that
$$[u_1,w_1]=0$$ and the associated symmetric Leibniz algebra is given by
$$
\left\{
\begin{array}{l}
u_1 \ast u_1=w_1, \\
u_1 \ast w_1= w_1 \ast u_1=0.
\end{array}
\right.$$

\medskip
\noindent{\bf Remark. Generalization in any dimension.} We can generalize this case, considering $\dim A=n, \dim A_\bullet ^2= 1.$ We consider a basis $\{u_1,v_1,\cdots,v_{n-2},w_1\}$ such that $\{v_1,\cdots,v_{n-2},w_1\}$ is a basis of the center of $(A,\bullet)$ and $\{w_1\}$ a basis of $A_\bullet^2$. The non trivial product of $(A,\bullet)$ is
$$u_1 \bullet u_1= w_1.$$
We define  the Lie algebra structures $(A,[\;,\,])$ on the basis $\{u_1,v_1,\cdots,v_{n-2},w_1\}$ as follow:
\begin{itemize}
\item We consider on the $\mathbb{K}$-vector space $\g=\mathbb{K}\{v_1,\cdots,v_{n-2}\}$ any $(n-2)$-dimensional Lie algebra structure with a nontrivial $H^2(\g,\g)$ where $H^*(\g,\g)$ is the Chevalley-Eilenberg cohomology of $\g$. 
\item Let $\theta$ be a nontrivial $2$-form on $\g$ (which exists because $H^2(\g,\g) \neq 0)$. We consider a central extension $\g_1=\g \oplus \K\{w_1\}$ of $\g$ associated with $\theta$, that is the bracket of $\g_1$ is defined from the bracket of $\g$ and by the relation 
$$[v_i,v_j]_{\g_1}=[v_i,v_j]_{\g}+\theta(v_i,v_j)w_1$$ for $i=1,\cdots,n-2$.
\item Let $g$ be a derivation of $\g_1$ satisfying $g(w_1)=0.$ 
 We consider an extension of $\g_1$ by the derivation $g$, that is if $A=\K\{u_1\}\oplus \g_1$, the bracket of $A$ is that of $\g_1$ and $[u_1,X]=g(X)$ for any $X \in \g_1$. 
\end{itemize}
Thus the Lie structure obtained on $A$ satisfies the required conditions.

\noindent{Examples.}
\begin{enumerate}
\item If $\dim \g=1$ that is $\g=\K\{v_1\}$, the Lie algebra $\g_1=\{v_1,w_1\}$ is abelian. Thus the Lie bracket on $(A,[\;,\,])$ is given by 
$$[u_1,v_1]=av_1+bw_1.$$
In this case we obtain the $3$-dimensional symmetric Leibniz algebras whose multiplication satisfies
$$
\left\{
\begin{array}{l}
u_1 \ast u_1= w_1,\\
 u_1 \ast v_1=-v_1 \ast u_1=av_1+bw_1,
\end{array}
\right.$$
and other products are equal to zero.
\item If $\dim \g=2$ and $\g=\K\{v_1,v_2\}$ is abelian. Then $\g_1=\K\{v_1,v_2,w_1\}$ is the Heisenberg algebra:
$$[v_1,v_2]=w_1$$
and considering a general derivation $g$ of $\g_1$ we obtain the Lie bracket of $(A,[\; , \,])$:
$$[u_1,v_1]=a_1v_1+b_1v_2+c_1w_1, \ [u_1,v_2]=a_2v_1+b_2v_2+c_2w_1, \ [v_1,v_2]=w_1.$$
We obtain the following $4$-dimensional symmetric Leibniz algebras
$$
\left\{
\begin{array}{l}
u_1 \ast u_1 =w_1,\\ u_1\ast v_1=-v_1\ast u_1=a_1v_1+b_1v_2+c_1w_1,\\ 
 u_1\ast v_2=-v_2\ast u_1=a_2v_1+b_2v_2+c_2w_1, \\
 v_1 \ast v_2=-v_2 \ast v_1 =w_1.
\end{array}
\right.$$

\noindent{\bf Particular case: the four-dimensional oscillator Lie algebra.} This case corresponds to $a_2=-b_1=1$ and the other parameters equal to zero. Then the Lie bracket of $(A,[\; , \,])$ is given by
$$[u_1,v_1]=v_2, \ [u_1,v_2]=-v_1, \ [v_1,v_2]=w_1.$$
This Lie algebra is usually called the oscillator Lie algebra. It is a linear Lie algebra whose elements are the matrices
$$
\left(
\begin{array}{cccc}
0 & -z & y & 2t	\\
0&0&-x & y\\
0 & x & 0 &z\\
0&0&0&0
\end{array}
\right)
$$
 The corresponding symmetric Leibniz algebra is given by
$$
\left\{
\begin{array}{l}u_1 \ast u_1 =w_1,\\ u_1\ast v_1=-v_1\ast u_1=-v_2,\\  u_1\ast v_2=-v_2\ast u_1=v_1, \\
 v_1 \ast v_2=-v_2 \ast v_1 =w_1.
\end{array}
\right.$$
\end{enumerate}
We find again the result proved in\cite{A.B} in which one shows that the oscillator Lie algebra can be endowed with a symmetric Leibniz algebra structure and with a Poisson algebra structure.
\end{enumerate}
\medskip

\noindent{\bf 2. $\dim A=3.$} 
\begin{enumerate}
\item Assume that $\dim A=3$ and the multiplication $\bullet$ is given by
$$u_1 \bullet u_1= u_2 \bullet u_2= w_1.$$
The corresponding Lie bracket satisfies $[u_i,w_1]=0$. We put $[u_1,u_2]=\alpha w_1$. We obtain the symmetric Leibniz algebra
$$
\left\{
\begin{array}{l}
u_1 \ast u_1= u_2\ast u_2=w_1\\
u_1 \ast u_2=-u_2 \ast u_1= \alpha w_1.$$
\end{array}
\right.
$$

\item
Assume that  the multiplication $\bullet$ is given by
$$u_1 \bullet u_2= u_2 \bullet u_1= w_1.$$
The corresponding Lie bracket satisfies $[u_i,w_1]=0$. We put $[u_1,u_2]=\alpha w_1$. We obtain the symmetric Leibniz algebra
$$
\left\{
\begin{array}{l}
u_1 \ast u_1= u_2\ast u_2=0\\
u_1 \ast u_2=(1+\alpha) w_1\\
u_2 \ast u_1= (1-\alpha) w_1.$$
\end{array}
\right.
$$
\item If we assume that $(A,\bullet)$ is decomposable with a non trivial product $\bullet$, then $\{u_1,v_1,w_1\}$ is the basis of $A$ and the corresponding symmetric Leibniz algebra is given by
$$
\left\{
\begin{array}{l}
u_1 \ast u_1= w_1,\\
 u_1 \ast v_1=-v_1 \ast u_1=av_1+bw_1.
\end{array}
\right.$$
We find the result given in the example (a) in the previous case.

\end{enumerate}

\medskip

\noindent{\bf 3. $\dim A=4.$} 
\begin{enumerate}
\item  $A=\K\{u_1,u_2,u_3 \} \oplus \K\{w_1\}$ and
$$u_1 \bullet u_1 = w_1, \ u_2 \bullet u_2 =w_1, \ u_3 \bullet u_3 = w_1.$$
In this case we have
$$[u_i,u_j]=\alpha_{i,j}w_1, \ \ 1 \leq i,j \leq 3.$$
The Lie algebra $(A,[\;,\,])$ is a nilpotent Lie algebra isomorphic to the direct sum of the $3$-dimensional Heisenberg algebra with a $1$-dimensional abelian Lie algebra. The corresponding symmetric Leibniz algebras are given by
$$
\left\{
\begin{array}{l}
u_1 \ast u_1= u_2 \ast u_2=u_3 \ast u_3=w_1,\\
 u_1 \ast u_2=-u_2 \ast u_1=\alpha_{1,2}w_1,\\
 u_1 \ast u_3=-u_3 \ast u_1=\alpha_{1,3}w_1,\\
 u_2 \ast u_3=-u_3 \ast u_2=\alpha_{2,3}w_1.\\
\end{array}
\right.$$

\item $A=\K\{u_1,u_2 \} \oplus \K\{w_1,w_2\}$ and
$$u_1 \bullet u_1 = w_1, \ \ u_1 \bullet u_2=u_2 \bullet u_1=w_2.$$
In this case, the non trivial Lie bracket is
$$[u_1,u_2]=\alpha w_1+\beta w_2$$
that is $(A,[\;,\,])$ is the Heisenberg algebra or the abelian Lie algebra. We deduce the corresponding symmetric Leibniz algebra:
$$
\left\{
\begin{array}{l}
u_1 \ast u_1= w_1,\\
 u_1 \ast u_2=\alpha w_1 + (\beta+1)w_2,\\
u_2 \ast u_1=-\alpha w_1 - (\beta-1)w_2.\\
\end{array}
\right.$$

\item $A=\K\{u_1,u_2 \} \oplus \K\{w_1,w_2\}$ and
$$u_1 \bullet u_1 = w_1, \ \ u_2 \bullet u_2 =  w_1, \ \ u_1 \bullet u_2=u_2 \bullet u_1=w_2$$
In this case also, we have
$$[u_1,u_2]=\alpha w_1+\beta w_2$$
and the corresponding symmetric Leibniz algebra is
$$
\left\{
\begin{array}{l}
u_1 \ast u_1=u_2 \ast u_2= w_1,\\
 u_1 \ast u_2=\alpha w_1 + (\beta+1)w_2,\\
u_2 \ast u_1=-\alpha w_1 - (\beta-1)w_2.\\
\end{array}
\right.$$

\item Now we consider that $(A,\bullet)$ is decomposable. The first case is when $(A,\bullet)$ is trivial that is $A=\K\{v_1,v_2,v_3,v_4\}$ that is $v_i \bullet v_j=0.$ In this case $(A,[\;,\,])$ is any $4$-dimensional Lie algebra and the symmetric Leibniz algebra coincides with this Lie algebra, the product $\ast$ is then skew-symmetric.

\item $A=\K\{u_1\} \oplus \K\{v_1,v_2\} \oplus \K\{w_1\}$. This is equivalent to write
$$u_1 \bullet u_1 = w_1.$$
In this case, we can  have 
$$\begin{array}{l}
\medskip
\lbrack u_1,v_j]=\ds\sum_{k=1}^2 E_{1,j}^kv_k+ F_{1,j}w_1, \ \ j=1,2\\
\medskip
\lbrack v_1,v_2]=\ds \sum_{k=1}^2 G_{1,2}^kv_k+ H_{1,2}w_1.\\
\end{array}
$$
This case has also been studied in a previous example. In particular, we find the $4$-dimensional oscillator Lie algebra.

\item $A=\K\{u_1,u_2\} \oplus \K\{v_1\} \oplus \K\{w_1\}$. This is equivalent to write
$$u_1 \bullet u_1 =u_2 \bullet u_2 = w_1.$$
In this case, the Lie bracket writes :
$$\begin{array}{l}
\medskip
\lbrack u_1,u_2]=\ds C_{1,2}v_1+ D_{1,2}w_1, \\
\medskip
\lbrack u_i,v_1]=\ds  E_{i,1}v_1+ F_{i,1}w_1, \ \ i=1,2\\
\end{array}
$$
with the Jacobi condition
$$F_{1,1}E_{2,1}-F_{2,1}E_{1,1}=0.$$
We deduce the following symmetric Leibniz algebras
$$
\left\{
\begin{array}{l}
u_1 \ast u_1=u_2 \ast u_2= w_1,\\
 u_1 \ast u_2=-u_2 \ast u_1= C_{1,2}v_1+ D_{1,2}w_1\\
u_1 \ast v_1=-v_1 \ast u_1=E_{1,1}v_1+ F_{1,1}w_1,\\
u_2 \ast v_1=-v_1 \ast u_2=E_{2,1}v_1+ F_{2,1}w_1.\\
\end{array}
\right.$$
\end{enumerate}

\section{Deformation quantization of Poisson algebras in a symmetric Leibniz formal deformation}

As it was recalled in the introduction, the formal deformations of commutative weakly associative algebras gives a construction of nonassociative Poisson algebras. An interesting case corresponds to a formal deformation of an associative commutative algebra in the category of weakly associative algebras which gives a classical Poisson algebra and enlarges the spectrum of deformation quantization. Moreover the class of weakly associative is the only one which permits to construct such a Poisson algebra.

In this paragraph we investigate the formal deformation of commutative symmetric Leibniz algebras in the class of symmetric Leibniz algebras.
 Let 
$(A, \ast)$ be  a commutative symmetric Leibniz   algebra. A formal deformation $(A[[t]], \ast_t)$ is given by a symmetric Leibniz formal product which can be represented by a formal series 
$$X \ast_t Y=X \ast Y+\sum_{i}t^i\varphi_i(X,Y)$$
for any $X,Y \in A,$ provides $A[[t]]$ with a symmetric Leibniz algebra structure. We then have 
$$
\left\{
\begin{array}{l}
X\ast_t(Y\ast_t Z) - (X\ast_t Y)\ast_t Z - Y\ast_t(X\ast_t Z)=0,\\
Y\ast_t (Z\ast_t X) - (Y\ast_t Z)\ast_t X + (Y\ast_t X)\ast_t Z =0
\end{array}
\right.
$$
This system implies 
\begin{itemize}
\item at the order 0 that $\ast$ is a symmetric Leibniz product,
\item at the order 1 that the linear bilinear map $\varphi_1$ satisfies
$$\left\{
\begin{array}{ll}
\delta \varphi_1^{(1)}
(X,Y,Z)&=\varphi_1(X, Y\ast Z)-\varphi_1(X\ast Y, Z)-\varphi_1(Y,X\ast Z)\\
&+ X\ast \varphi_1(Y,Z)-\varphi_1(X,Y)\ast Z -Y \ast \varphi_1(X,Z)=0 \\
\delta\varphi_1^{(2)}
(X,Y,Z)&=\varphi_1(Y, Z\ast X)-\varphi_1(Y\ast Z, X)+\varphi_1(Y\ast X, Z)\\
& +Y\ast \varphi_1(Z,X)-\varphi_1(Y,Z)\ast X + \varphi_1(Y,X)\ast Z=0
\end{array}
\right.$$
\item at the order 2 that the linear bilinear map $\varphi_2$ satisfies
$$\left\{
\begin{array}{ll}
\mathcal{A}_{\varphi_1}
(X,Y,Z)-\varphi_1( Y,\varphi_1(X, Z))+\delta\varphi_2^{(1)} 
(X,Y,Z)=0\\
\mathcal{A}_{\varphi_1}
(Y,Z,X)+\varphi_1( \varphi_1(Y,X), Z)+\delta \varphi_2^{(2)}
(X,Y,Z)=0
\end{array}
\right.
$$
which is equivalent, since $\ast$ is supposed to be commutative,  to 
$$\left\{
\begin{array}{ll}
\mathcal{A}_{\varphi_1}
(X,Y,Z)+\mathcal{A}_{\varphi_1}(Y,Z,X)-\mathcal{A}_{\varphi_1}(Y,X,Z)+
\delta \varphi_2^{(1)}(X,Y,Z)+\delta \varphi_2^{(2)}
(X,Y,Z)=0\\
\varphi_1(\varphi_1(X, Z),Y)+\varphi_1( Y,\varphi_1(X, Z))+\delta \varphi_2^{(2)}(Z,X,Y)-\delta \varphi_2^{(1)} (X,Y,Z)=0.
\end{array}
\right.
$$
If we put $\delta \varphi_2=\delta \varphi_2^{(1)}+\delta \varphi_2^{(2)}$ then
$$\left\{
\begin{array}{rl}
\delta \varphi_2
(X,Y,Z)=& \varphi_2(X,Y\ast Z)-\varphi_2(X,Y)\ast Z-
Y \ast \varphi_2(X,Z)-\varphi_2
(Y\ast Z,X)\\
& +Y
\varphi_2( Z,X)+\varphi_2( Y,X) \ast Z\\
= &2 \psi_{\varphi_2}(X,Y\ast Z)-2 \psi_{\varphi_2}(X,Y) \ast Z-2Y\ast  \psi_{\varphi_2}(X,Z)
\end{array}
\right.
$$
where $2\psi_{\varphi_2}(U,V)=\varphi_2(U,V)-\varphi_2(V,U)$ for $U,V \in A$.
Then $$\delta \varphi_2(X,Y,Z)-\delta \varphi_2(X,Z,Y)=0.$$
We deduce that
$$\begin{array}{c}
\mathcal{A}_{\varphi_1}
(X,Y,Z)+\mathcal{A}_{\varphi_1}
(Y,Z,X)-\mathcal{A}_{\varphi_1}
(Y,X,Z)-
\mathcal{A}_{\varphi_1}
(X,Z,Y)\\-\mathcal{A}_{\varphi_1}
(Z,Y,X)+\mathcal{A}_{\varphi_1}
(Z,X,Y)=0
\end{array}$$
which is equivalent to say that $\varphi_1$ is Lie-admissible or equivalently that $\psi_{\varphi_1}$ is a Lie bracket.
Since $\delta \varphi_1^{(1)}=\delta \varphi_1^{(2)}=0$, from this calculus we deduce also  that
$$
\begin{array}{rl}
\ds \frac{ \delta \varphi_1
(X,Y,Z)}{2}= & \psi_{\varphi_1}(X,Y\ast Z)- \psi_{\varphi_1}(X,Y) \ast Z-Y\ast  \psi_{\varphi_1}(X,Z)=0
\end{array}
$$
then $(A,\ast, \psi_{\varphi_1})$ is a nonassociative Poisson algebra defined by the commutative symmetric Leibniz algebra $(A,\ast).$ So we find the same result concerning weakly associative algebra in the symmetric Leibniz context. 
\end{itemize}  
The Leibniz identity between $\psi_{\varphi_1}$ and $\ast$ is equivalent to $\delta \varphi_1^{(1)}+\delta \varphi_1^{(2)}=0$. Since we have  $\delta \varphi_1^{(1)}=\delta \varphi_1^{(2)}=0$ we have to look the consequences of $\tilde{\delta} \varphi_1=\delta \varphi_1^{(1)}-\delta \varphi_1^{(2)}=0$. We have
$$
\begin{array}{ll}
\tilde{\delta} \varphi_1
(X,Y,Z)&=\varphi_1(X, Y\ast Z)-\varphi_1(X\ast Y, Z)-\varphi_1(Y,X\ast Z)\\
&+ X\ast \varphi_1(Y,Z)-\varphi_1(X,Y)\ast Z -Y \ast \varphi_1(X,Z) \\
&-\varphi_1(Y, Z\ast X)+\varphi_1(Y\ast Z, X)-\varphi_1(Y\ast X, Z)\\
& -Y\ast \varphi_1(Z,X)+\varphi_1(Y,Z)\ast X - \varphi_1(Y,X)\ast Z\\
&=2\rho_{\varphi_1}(X, Y\ast Z)-2\varphi_1(X\ast Y, Z)-2\varphi_1(Y,X\ast Z)\\
&+2 X\ast \varphi_1(Y,Z)-2\rho_{\varphi_1}(Y,X)\ast Z-2Y\ast\rho_{\varphi_1}(Z,X)
\end{array}
$$
Since $\varphi_1=\rho_{\varphi_1}+\psi_{\varphi_1}$, we obtain
$$
\begin{array}{ll}
\tilde{\delta} \varphi_1
(X,Y,Z)&=2\rho_{\varphi_1}(X, Y\ast Z)-2\rho_{\varphi_1}(Y,X)\ast Z-2Y\ast\rho_{\varphi_1}(Z,X)-2\rho_{\varphi_1}(X\ast Y, Z)\\
&-2\rho_{\varphi_1}(Y,X\ast Z)+
2 X\ast \rho_{\varphi_1}(Y,Z)-2\psi_{\varphi_1}(X\ast Y, Z)
-2\psi_{\varphi_1}(Y,X\ast Z)\\
&+2 X\ast \psi_{\varphi_1}(Y,Z)
\end{array}
$$
To simplify the writing, because any confusion is possible, we write $XY$ instead $X \ast Y$ and also $\rho,\psi$ in place of $\rho_{\varphi_1},\psi_{\varphi_1}$. The previous identity writes
\begin{eqnarray}\label{rho}
0=&\rho(X, YZ)-\rho(Y,X)Z-Y\rho(Z,X)-\rho(X Y, Z)
-\rho(Y,XZ)+
 X \rho(Y,Z)-\psi(XY, Z) \nonumber\\
&-\psi(Y,X Z)
+ X\psi(Y,Z) 
\end{eqnarray}

Permuting $X$ and $Z$ we obtain
$$\begin{array}{ll}
0=&\rho(Z, YX)-\rho(Y,Z)X-Y\rho(X,Z)-\rho(Z Y, X)
-\rho(Y,ZX)+
 Z \rho(Y,X)-\psi(ZY, X)\\
&-\psi(Y,ZX)
+ Z\psi(Y,X)
\end{array}
$$
Adding these two relations
$$\begin{array}{ll}
0=&-2Y\rho(X,Z)
-2\rho(Y,ZX)-\psi(ZY, X)-2\psi(Y,ZX)
+ Z\psi(Y,X)-\psi(XY, Z)\\&+ X\psi(Y,Z)
\end{array}
$$

 We have seen that $\psi_{\varphi_1}$ is a Lie bracket satisfying the Leibniz identity with respect the associative commutative multiplication $\bullet$ but also with respect the multiplication $\ast$. We deduce that
 
$$-\psi(ZY, X)-2\psi(Y,ZX)
+ Z\psi(Y,X)-\psi(XY, Z)+ X\psi(Y,Z)=-2\psi(Y,XZ).$$
So we obtain
$$Y\rho(X,Z)
+\rho(Y,ZX)+\psi(Y,XZ)=0.$$

\medskip

The identity (\ref{rho}) becomes
$$
-\psi(X,YZ)
+\psi(Z,XY)
+\psi(Y,XZ)
-\psi(XY, Z)-\psi(Y,X Z)
+ X\psi(Y,Z)$$
that is 
$$
-\psi(X,YZ)
+2\psi(Z,XY)
+ X\psi(Y,Z)=0.
$$
From the Leibniz identity between $\psi$ and $\ast$, this is equivalent to
$$\psi(Y,Z)X+\psi(X,Y)Z-3\psi(Z,X)Y=0.$$
Then we have also
$$\psi(Z,Y)X+\psi(X,Z)Y-3\psi(Y,X)Z=0.$$
If we add these two last identities we obtain
$$4\psi(X,Y)Z-4\psi(Z,X)Y=0.$$

\begin{theorem}
Let $(A,\ast)$ be a  commutative symmetric Leibniz algebra and $\ast_t=\ast +\sum _{i\geq 1} t^i \varphi_i$ a symmetric Leibniz formal deformation of $\ast$ (that is $(A[[t]],\ast_t)$ is a symmetric Leibniz algebra). Then if $\psi_{\varphi_1}$ is the skew-symmetric map attached with $\varphi_1$, the algebra $(A,\ast, \psi_{\varphi_1})$ is a nonassociative Poisson algebra and the Lie bracket $\psi_{\varphi_1}$ satisfies
$$\psi_{\varphi_1}(X,Y)\ast Z-\psi_{\varphi_1}(Z,X)\ast Y=0$$
for any $X,Y,Z \in A$.
\end{theorem}

\medskip

\noindent Remark. A commutative associative algebra is always weakly associative. But it is a symmetric Leibniz algebra if and only if it is $2$-step nilpotent. Then the previous theorem permits to construct from a commutative associative $2$-step nilpotent algebra a classical Poisson algebra  via a formal symmetric Leibniz deformation process.

\end{document}